\documentclass[fleqn,10pt]{wlscirep}
\usepackage[utf8]{inputenc}
\usepackage[T1]{fontenc}

\usepackage{algorithm}
\usepackage{algorithmic}
\usepackage{subcaption}
\usepackage{float}
\usepackage{lineno}
\usepackage{diagbox}
\usepackage{amsmath, amssymb, amsbsy, amsthm, epigraph, longtable, soul, hyperref}
\usepackage{graphicx, xcolor, url, tikz}
\usetikzlibrary{fit}

\newcommand{\R}{\mathbb{R}}

\renewcommand{\phi}{\varphi}
\renewcommand{\theta}{\vartheta}

\usepackage{etoolbox}
\newcommand*\linenomathpatch[1]{%
  \cspreto{#1}{\linenomath}%
  \cspreto{#1*}{\linenomath}%
  \csappto{end#1}{\endlinenomath}%
  \csappto{end#1*}{\endlinenomath}%
}

\linenomathpatch{equation}
\linenomathpatch{gather}
\linenomathpatch{multline}
\linenomathpatch{align}
\linenomathpatch{alignat}
\linenomathpatch{flalign}
\newcommand{\lpm}{\begin{pmatrix}}
\newcommand{\rpm}{\end{pmatrix}}

\theoremstyle{plain}

\theoremstyle{definition}

\theoremstyle{remark}
   
\begin{document}
\title{Modelling the Spatial Spread of COVID-19 in a German District using a Diffusion Model}
\author[1,*]{Moritz Schäfer}
\author[1]{Peter Heidrich}
\author[1]{Thomas Götz}



\affil[1]{Mathematical Institute, University of Koblenz, 56070 Koblenz, Germany}


\begin{abstract}

This study focuses on modeling the local spread of COVID-19 infections. As the pandemic continues and new variants or future pandemics can emerge, modelling the early stages of infection spread becomes crucial, especially as limited medical data might be available initially. Therefore, our aim is to gain a better understanding of the diffusion dynamics on smaller scales using partial differential equation (PDE) models.

The article focuses on a single German district, \textit{Birkenfeld} in Rhineland-Palatinate, during the second wave of the pandemic in autumn 2020 and winter 2020--21. This district is characterized by its (mainly) rural nature and daily commuter movements towards metropolitan areas. Medical data from both the Robert-Koch-Institute (RKI) and the Birkenfeld district government for parameter estimation is used.

A basic reaction-diffusion model used for spatial COVID spread, which includes compartments for susceptibles, exposed, infected, recovered, and the total population, is used to describe the spatio-temporal spread of infections. The transmission rate, recovery rate, initial infected values, detection rate, and diffusivity rate are considered as parameters to be estimated using the reported daily data and least square fit. The optimization is performed using two different methods, the Metropolis algorithm and the adjoint method.


\end{abstract}

\newpage

\maketitle

\section{Introduction}

At the beginning of January 2020, the COVID-19 virus began to spread throughout mainland China, with the consequences that we all have experienced in the last three years. Initially, the number of cases was limited to single clusters in a limited number of locations, but later on expanded throughout the country. In previous studies, we have investigated the macroscopic impact of the epidemic using an \textit{SIR}-model for all cases in Germany. For this, we have used classical differential models such as the \textit{SEIRD}- (susceptible-exposed-infected-recovered-dead) models to describe the spread of infections during the first wave (cf. Heidrich et al. \cite{HeS20}), as well as the impact of travelers on disease dynamics in summer 2020 (cf. Schäfer et al. \cite{Sch22}), both with a strong emphasis on parameter estimation. 

In this study, we aim to model the local spread of infections using PDE (partial differential equation) models to gain a better understanding of the diffusion on smaller scales, similar to the work of Viguerie et al. \cite{Vig21} and Wang and Yamamoto \cite{Wan20}. Kuehn and Mölter \cite{Kue22} also investigated the impact of transportation on epidemics using two coupled models, a static epidemic network and a dynamic transportation network, with non-local, fractional transport dynamics. Logeshwari et al. \cite{Log20} also provide a fractional PDE model of the spatial spread of COVID-19. A general challenge with diffusion-based PDE models is that diffusion of all compartments leads to unwanted diffusion in the total population, which we aim to avoid. This study uses data down to the municipality level for parameter estimation of the model.

\begin{figure}[H]
\begin{center}
		\includegraphics[width=.75\linewidth]{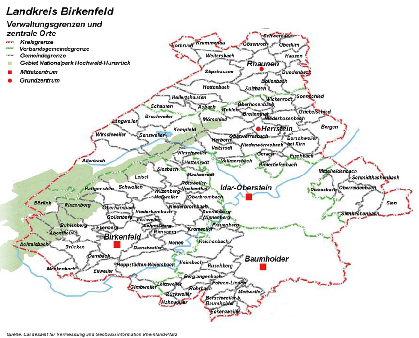}
		\centering
	\caption{Map of Birkenfeld County with administrative boundaries of central localities \cite{Lan}.}
	\label{fig:data1A}
	\end{center}
\end{figure}

We focus our numerical problem on a single district, the district of Birkenfeld in southwestern Germany within the state of Rhineland-Palatinate. Approximately 81,000 people live there in an area of about 780 km$^2$. The largest city within the district is Idar-Oberstein, with about 28,000 inhabitants in an area of about 92 km$^2$. The remaining people live in the municipalities of Birkenfeld, Baumholder, and Herrstein-Rhaunen. Within a 1.5 hour drive via federal highways and freeways, the following metropolitan areas can be reached: Mainz, Trier, Koblenz, Kaiserslautern, Saarbrücken, and Frankfurt. In addition, the Frankfurt-Hahn airport is located in the neighboring Rhein-Hunsrück district to the northeast. The region is very rural, with daily commuter movements common in the direction of the aforementioned metropolitan areas. The region is also visited by tourists due to the gemstone industry and trade in Idar-Oberstein, the numerous hiking routes, and the nearby Hunsrück-Hochwald National Park. A map of this district can be found in Fig. \ref{fig:data1A}.

\begin{figure}[H]
	\begin{subfigure}{.5\textwidth}
		\centering
		\includegraphics[width=1\linewidth]{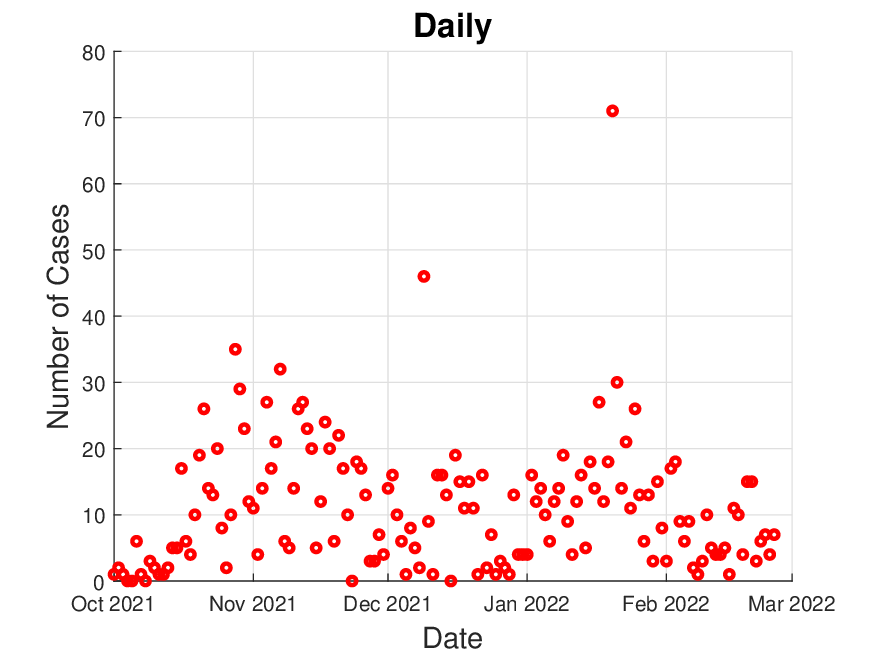}
		\end{subfigure}%
	\begin{subfigure}{.5\textwidth}
		\centering
		\includegraphics[width=1\linewidth]{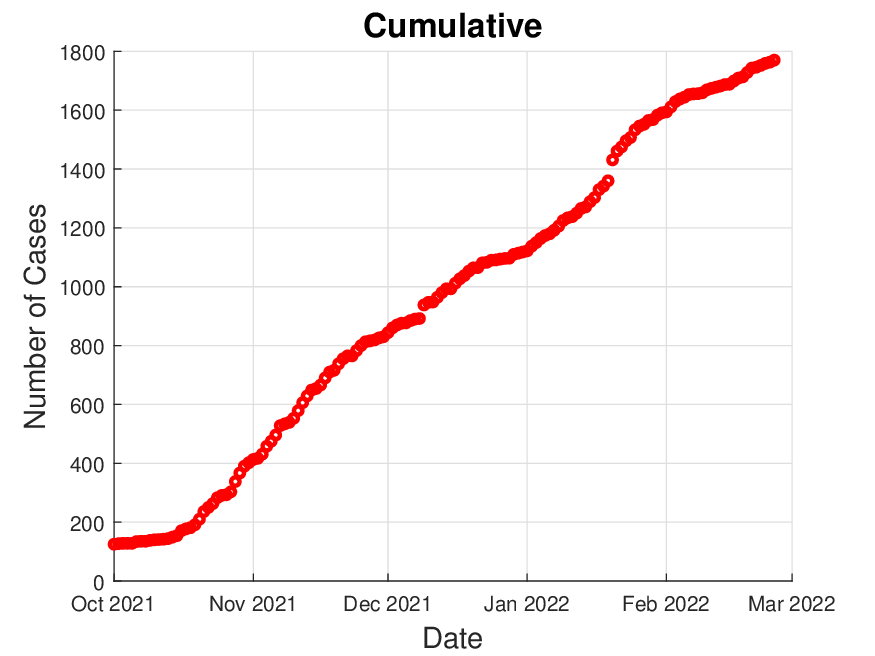}
		\end{subfigure}
	\caption{Confirmed new daily cases (left) and cumulative confirmed cases (right) with COVID-19 in Birkenfeld from October 1, 2020 until February 25, 2021 according to RKI \cite{RKI2}.  }
	\label{fig:data1C}
\end{figure}

The first COVID-19 case in the district of Birkenfeld was confirmed on March 16, 2020. Until June 30, 2020, there were only 90 registered cases in the entire district. However, the number of cases increased during the second wave in autumn/winter 2020/21, with a cumulative 2,513 cases confirmed until March 31, 2021 \cite{RKI2}. As of October 2022, over 32,000 cases have been counted in the district. For several reasons, we restrict our research to the data from the second wave: there was only a limited number of previous infections, but comparatively high infection numbers further on; also, no vaccines were available until the beginning of 2021, as well as a very small amount of persons having been exposed to the disease multiple times. Finally, the political lockdown restrictions, particularly in November and December 2020, led to a slower mixing of cases between different districts. This inter-district mobility is not taken into account by our models. We consider daily infection data from the district and all of its municipalities in the time frame from October 1, 2020 to February 25, 2021. the cumulative number of infections in the Birkenfeld district is depicted in Fig. \ref{fig:data1C}. In this figure, the first step visible in the data between end of December and beginning of January is mainly caused by a delay of registration of case numbers due to Christmas holidays. The reasons for second one in mid January are unknown and could be related to registration delays within the district. An important fact is that there are only one detected initial infection case on October 1, situated in the city of Idar-Oberstein, so there were no detected cases in the rest of the entire district.

As the pandemic continues to spread globally (as of autumn 2022), the possibility of new variants of the virus or future pandemics highlights the importance of modeling the spread of infections, particularly in their early stages when limited medical data is available. The accuracy of these models depends greatly on the used parameters and the corresponding data. In this study, we present \textit{SEIRD}-models, which are commonly used in epidemiological simulations, and estimate their parameters using data from the Robert-Koch-Institute (RKI) \cite{RKI2} and private communications with the Birkenfeld district government\cite{priv}. We perform the estimations using both adjoint and Metropolis methods and base them on a least square fit between the model output and the reported daily data.

\section{Materials and methods}

\subsection{PDE models}

To model the spatial COVID spread in the presented areas, we use an epidemiological reaction-diffusion model. For this purpose, we consider a corresponding spatial area $(x,y) \in \Omega$ in a time period $t \in \mathcal T := \left[ 0 , t_{end} \right]$. We are looking for a function $u : V \rightarrow \mathbb{R}^m$ with $V=\Omega \times \mathcal T $, which is twice continuously differentiable on $\Omega$ and once continuously differentiable on $\mathcal T $, briefly $u \in \mathcal{C}(V,\mathbb{R}^m)^{2,1}$. The following PDE system has to be fulfilled:
\begin{subequations}
\label{eq:u_PDE_}
    \begin{alignat}{2}
        \partial_t u    &= \kappa \Delta_{x,y} u + f(u) \, , \qquad &&
        \\
        u &= u_0 \, ,   &&  t = 0 \, ,
        \\
        \partial_\nu u  &= 0 \, ,   &&(x,y) \in \partial \Omega \, .
    \end{alignat}
\end{subequations}
Here, $\partial_t u$ stands for the component--wise derivative of $u$ in the direction of time, i.e., $\partial_t u = \left( {\partial_t u_1} ,...,{\partial_t u_m} \right)^T$ and $\Delta_{x,y} = \left({\partial_{xx} u_1} + {\partial_{yy} u_1},..., {\partial_{xx} u_m} + {\partial_{yy} u_m} \right)^T$ for the Laplace operator in $\Omega$. The parameter $\kappa$ describes the diffusivity of the system and the function $f(u)$ contains the epidemiological component(s). As an initial condition, at time $t=0$ a function $u_0 : \Omega \rightarrow \mathbb{R}^m$ is used with $u(x,y,t=0)=u_0(x,y)$. In addition, Neumann boundary conditions are used, where $\partial_\nu u = \left( {\partial_\nu u_1},..., {\partial_\nu u_m} \right)^T$ stands for the derivative in the direction of the outward pointing unit normal $\nu$ and $\partial \Omega$ stands for the boundary of $\Omega$. In terms of context, the latter means that no individual can leave or enter the territory $\Omega$. This seems strange at first, since the district of Birkenfeld in practice can be left or entered by land. On the other hand, we are looking at data sets from a period when profound containment measures had already been taken in the region and social measures, including a significant reduction inter-district mobility, were already implemented.

For epidemiological modelling, we make use of a variant of the \textit{SIR} model introduced by Kermack and McKendrick~\cite{KMK91}, the \textit{SEIR}-model, and consider compartments as functions $S,E,I,R,N \in \mathcal{C}(V,\mathbb{R})^{2,1}$, which have the following meanings:
\begin{itemize}
    \item Susceptibles $S$: Depending on the transmission route, these individuals can become infected with the infectious disease when contact occurs. 
    \item Exposed $E$: The corresponding indiviuals have already ingested the pathogen, but are not yet infectious because they are still in the latency period.
    \item Infected $I$: These individuals are infected with the disease and infectious. Contact with a susceptible individual can therefore lead to transmission of the disease.
    \item Recovered $R$: After surviving an infection, individuals are considered recovered. These individuals can no longer transmit the disease or get infected. 
    \item Population $N$: The total number of individuals.
\end{itemize}
For instance, $I(x,y,t)$ indicates the number of infected individuals in the spatial coordinate $(x,y) \in \Omega$ at time $t \in \mathcal T $. Based on these presented groups, different epidemiological models can now be derived. We present here the PDE systems of three common models in Tab. \ref{tab:models_}.
\begin{table}[h!btp]
\caption{Basic examples of epidemiological compartment models with flow chart and PDE system.}
\centering
	\begin{tabular}{c|c|c}
		\hline
		Model	&	Structure	&	PDE System
		\\
		\hline	
				&				&
		\\		
		\textit{SIS}	& 			\begin{tikzpicture}[node distance=2cm, scale=0.33\textwidth, baseline=0]
							\tikzstyle{comp} = [rectangle, rounded corners, minimum width=.7cm, minimum height=.7cm, text centered, draw=black]
							\tikzstyle{arrow} = [thick,->,>=stealth];
							\node (S) [comp, fill=green!30] {$S$};
							\node (I) [comp, right of=S, fill=red!60] {$I$};
							\draw[arrow] (S)-- node[above, pos=0.5]{\footnotesize{$\dfrac{\beta(t)}{N}\, SI$}} (I);
							\draw[arrow] (I)|- ++(0,-0.005) -| node[below, pos=0.25]{\footnotesize{$\gamma\, I$}} (S);
  							\end{tikzpicture}  							
  				&					{$\begin{aligned}
										\partial_t S & = \kappa_S \Delta_{x,y} S -\frac{\beta(t)}{N}SI + \gamma I, & S(t=0)=S_0,
										\\
										\partial_t I & = \kappa_I \Delta_{x,y} I + \frac{\beta(t)}{N}SI - \gamma I, & I(t=0)=I_0, 
										\\
										\partial_\nu S &= \partial_\nu I = 0,   &(x,y) \in \partial \Omega,
										\\
										N &= S + I\,.
									\end{aligned}$}	
		\\
		\hline	
				&				&
		\\		
		\textit{SIR}	&	\begin{tikzpicture}[node distance=2cm, scale=0.3\textwidth, baseline=0]
							\tikzstyle{comp} = [rectangle, rounded corners, minimum width=.7cm, minimum height=.7cm, text centered, draw=black]
							\tikzstyle{arrow} = [thick,->,>=stealth];
							\node (S) [comp, fill=green!30] {$S$};
							\node (I) [comp, right of=S, fill=red!60] {$I$};
							\node (R) [comp, right of=I, fill=green!70] {$R$};
 							\draw[arrow] (S)-- node[above, pos=0.5]{\footnotesize{$\dfrac{\beta(t)}{N}\, SI$}} (I);
  							\draw[arrow] (I) -- node[above, pos=0.5]{\footnotesize{$\gamma\, I$}} (R);
							\end{tikzpicture}			
							
						&	{$\begin{aligned}
										\partial_t S & = \kappa_S \Delta_{x,y} S -\frac{\beta(t)}{N}SI, & S(t=0)=S_0,
										\\
										\partial_t I & = \kappa_I \Delta_{x,y} I + \frac{\beta(t)}{N}SI - \gamma I, & I(t=0)=I_0, 
										\\
										\partial_t R & = \kappa_R \Delta_{x,y} R + \gamma I, & R(t=0)=R_0,
										\\
										\partial_\nu S &= \partial_\nu I = \partial_\nu R = 0,   &(x,y) \in \partial \Omega,
										\\
										N &= S + I + R\,.
									\end{aligned}$}	
		\\
		\hline		
				&				&	
		\\
		\textit{SEIR}	&			\begin{tikzpicture}[node distance=2cm, scale=0.3\textwidth, baseline=0]
							\tikzstyle{comp} = [rectangle, rounded corners, minimum width=.7cm, minimum height=.7cm, text centered, draw=black]
							\tikzstyle{arrow} = [thick,->,>=stealth];
							\node (S) [comp, fill=green!30] {$S$};
							\node (E) [comp, right of=S, fill=orange!30] {$E$};
							\node (I) [comp, right of=E, fill=red!60] {$I$};
							\node (R) [comp, right of=I, fill=green!70] {$R$};
 							\draw[arrow] (S)-- node[above, pos=0.5]{\footnotesize{$\dfrac{\beta(t)}{N}\, SI$}} (E);
 							\draw[arrow] (E) -- node[above, pos=0.5]{\footnotesize{$\theta\, E$}} (I);
  							\draw[arrow] (I) -- node[above, pos=0.5]{\footnotesize{$\gamma\, I$}} (R);
							\end{tikzpicture}				
							
						&	{$\begin{aligned}
										\partial_t S & = \kappa_S \Delta_{x,y} S -\frac{\beta(t)}{N}SI,  & S(t=0)=S_0,
										\\
										\partial_t E & = \kappa_E \Delta_{x,y} E + \frac{\beta(t)}{N}SI - \theta E,  & E(t=0)=E_0,
										\\
										\partial_t I & = \kappa_I \Delta_{x,y} I + \theta E - \gamma I,  & I(t=0)=I_0,
										\\
										\partial_t R & = \kappa_R \Delta_{x,y} R + \gamma I, & R(t=0)=R_0,
										\\
										\partial_\nu S &= \partial_\nu E = \partial_\nu I = \partial_\nu R = 0,   &(x,y) \in \partial \Omega,
										\\
										N &= S + E + I + R \, .
									\end{aligned}$}	
		\\
				&				&
		\\
		\hline
	\end{tabular}
\label{tab:models_}
\end{table} 
The derivation and precise functioning of spatial epidemiological models will not be explained in detail here; for this purpose, reference is made to e.g. Martcheva \cite{Mar15}. At the core of every epidemiological model is the so--called incidence term $\frac{\beta(t)}{N}SI$, which indicates how many individuals are newly infected with the disease in coordinate $(x,y)$ at time $t$. The incidence term depends on a time dependent transmission rate $\beta : [0,t_{end}] \rightarrow \R^+$. In simple models, this can also be assumed to be a constant parameter, but we assume that the transmission rate may fluctuate over the observed periods due to the stepwise restrictions on the population. The value of the transmission rate $\beta$ is generally unknown and must be adjusted using the data sets. Another parameter in the models is the recovery rate $\gamma$. This is the reciprocal of the time required on average for an individual to recover from the disease. Thus, if we assume that $t$ is in days and an individual takes $10$ days to recover, it holds $\gamma = \frac{1}{10}$. In addition, the \textit{SEIR}-model contains the parameter $\theta$, which is the reciprocal of the latency period, i.e. the time between the uptake of the pathogen into the body and the onset of infectiousness. For example, assuming three days, it holds $\theta = \frac{1}{3}$. It should be noted here that the latency period does not have to be congruent with the incubation period, as the latter indicates the period of time until the onset of the first symptoms. With regard to COVID--19 in particular, it has been shown that infectivity sets in even before the onset of symptoms (cf. He et al. \cite{HeL20}). Due to simplicity, we also chose $\kappa=\kappa_S=\kappa_E=\kappa_I=\kappa_R$.

We use the \textit{SEIR}-model as an example to show how the models are prepared for later data fitting. In the first step, we substitute $R = N-S-E-I$ and thus reduce the system to an \textit{SEI}-model. It should be noted, however, that also for $N$ a PDE has to be solved for which holds 
\begin{subequations}
    \begin{alignat}{2}
        \partial_t N &= \kappa \Delta_{x,y} N \, , \qquad &&
        \\
        N &= S_0 + E_0 + I_0 + R_0 \, , \qquad \qquad &&t = 0 \, ,
        \\
        \partial_\nu N &= 0 \, , &&(x,y) \in \partial \Omega \, .
    \end{alignat}
\end{subequations}
For this reason, we normalize the reduced \textit{SEI}-model by dividing all rows by $N$, assuming that the population density mathematically fulfills $N(x,y,t) > 0$ on $V$. Defining  $u_1:=\frac{S}{N}$, $u_2:=\frac{E}{N}$, $u_3:=\frac{I}{N}$ and $u:=(u_1,u_2,u_3)$, we obtain a system as in \eqref{eq:u_PDE_} with 
\begin{align}
\label{eq:f(u)_}
    f : \mathcal{C}(V,\mathbb{R}^3)^{2,1} \rightarrow \mathbb{R}^3, \qquad f(u)=\begin{pmatrix} -\beta(t)u_1u_3 \\ \beta(t)u_1u_3-\theta u_2 \\ \theta u_2 - \gamma u_3 \end{pmatrix} \, .
\end{align}
Tab. \ref{tab:f(u)_} summarizes the results for the presented models. 
\begin{table}[H]
    \centering
        \caption{Summary of $f(u)$ for the reduced and normalized models.}
    \label{tab:f(u)_}
    \begin{tabular}{c|cl}
        \hline
       
        \textit{SIS} $\rightarrow$ \textit{I }        &   $f : \mathcal{C}(V,\mathbb{R})^{2,1} \rightarrow \mathbb{R}$, 
                    &   $f(u)= \beta(t)(1-u)u - \gamma u$
        \\
        \hline
        \textit{SIR} $\rightarrow$ \textit{SI }        &   $f : \mathcal{C}(V,\mathbb{R}^2)^{2,1} \rightarrow \mathbb{R}^2$,
                    &   $f(u)=\begin{pmatrix}-\beta(t)u_1u_2 \\\beta(t)u_1u_2-\gamma u_2 \end{pmatrix}$
        \\
        \hline
        \textit{SEIR} $\rightarrow$ \textit{SEI}       &   $f : \mathcal{C}(V,\mathbb{R}^3)^{2,1} \rightarrow \mathbb{R}^3$,
                    &   $f(u)=\begin{pmatrix} -\beta(t)u_1u_3 \\ \beta(t)u_1u_3-\theta u_2 \\ \theta u_2 -       \gamma u_3 \end{pmatrix}$
        \\
        \hline
    \end{tabular}
\end{table}
To meaningfully include the biological context, it must hold $\gamma,\theta > 0$ and $\kappa \geq 0$ and, as initial condition, $u_0 \geq0$ in $\Omega$. In addition, we assume that there are infected individuals in the area $\Omega$, i.e. $\int_\Omega I_0 \, d\omega > 0$. For the reduced and normalized \textit{SEI}-model, it must then hold $\int_\Omega u_{3,0} \, d\omega > 0$, using the notation $u_j(x,y,0):=u_{j,0}$. The notation $\mathcal{N}_0 := \int_\Omega N(x,y,0) \, d \omega$ represents the total number of individuals at time $t=0$ in $\Omega$. It must be valid that $\mathcal{N}_0>0$. Moreover, we define the total population in the area $\Omega$ at time $t$ as $\mathcal{N}: [0,t_{end}] \rightarrow (0,+\infty)$ with $\mathcal{N}(t) := \int_\Omega N(x,y,t) \, d\omega$. Due to the Neumann boundary conditions, we receive using Gauss's \begin{align}
    \partial_t \mathcal{N}=\int_\Omega \partial_t N \, d \omega = \int_\Omega \kappa \,\Delta_{x,y} N \, d \omega = \int_{\partial \Omega} \kappa \,\partial_\nu N \, ds = 0.
\end{align} 
Thus, the total population in the domain $\Omega$ is constant with respect to time. Analytically, there exists a unique solution for each of the PDE systems \eqref{eq:u_PDE_} with the presented $f(u)$ in Table \ref{tab:f(u)_} in conjunction with the mentioned preconditions \cite{Bri86}.

Due to the formulation using diffusion for the total population, an equilibrium will only set in when the population density in the entire district is equal. Then, the temporal equilibrium will be analogous to the equilibrium of the system without diffusion.

As already mentioned, certain parameters of the model are known, such as $\gamma$ and $\theta$. The transmission rates $\beta_j$ and the diffusivity $\kappa$ are usually unknown. For the transmission rate, we assume that the time-dependancy is piecewise constant. Due to 'light' lockdown restrictions from November 2, 2020, and 'stricter' restrictions from December 17, 2020 to the end of the observed time interval, we assume three different time intervals as follows:
\begin{equation}
\label{eq:beta}
    \beta(t)= 
    \begin{cases}
    \beta_0 \, ,        &       0 \leq t < t_0 \, ,
    \\
    \beta_1 \, ,        &       t_0 \leq t < t_1 \, ,
    \\
    \beta_2 \, ,        &       t_1 \leq t \leq t_{end} \, .
    \end{cases}
\end{equation}In addition, due to noisy data sets, the initial conditions $u_0$ must also be adjusted. For that, we present two approaches to solve this in the following sections.

\subsection{Crank-Nicholson method for the SEIR-model}

For this purpose, we discretize the studied region $\Omega$ in $x$ and $y$ directions in equal equidistant step sizes $h_x$ and $h_y$, respectively. Also, the time interval $\mathcal T=[0,t_{end}]$ is divided into equidistant steps of length $\tau$. In the following we use the notation
$u_{i,j}^n = u(x_j,y_i,t_n)$. The Laplace operator is expressed by finite differences 
\begin{equation}
	\Delta u_{i,j}^n = \frac{1}{h_x^2} \left[ u_{i,j-1}^n - 2 u_{i,j}^n + u_{i,j+1}^n \right]
	                 + \frac{1}{h_y^2} \left[ u_{i-1,j}^n - 2 u_{i,j}^n + u_{i+1,j}^n \right]   
\end{equation}
and the Crank-Nicolson scheme reads as
\begin{equation}
\label{eq:cranknic}
    \frac{u_{i,j}^{n+1}-u_{i,j}^n}{\tau}=\frac 1 2 \left[ \kappa \Delta u_{i,j}^{n+1} + \kappa \Delta u_{i,j}^n + f(u_{i,j}^{n+1}) + f(u_{i,j}^n)  \right] \, .
\end{equation}
The goal is to transform this approach so that $u_{i,j}^{n+1}$ can be solved with a linear system of equations. The non-linearity of $f(u)$ is solved by evolving $f$ using Taylor expansion for small values $\tau$ around the current iteration value $u_{i,j}^n$ with $f(u_{i,j}^{n+1})=f(u_{i,j}(t^n + \tau))$:
 \begin{align}\label{eq:taylorf}
  \notag  f(u_{i,j}(t^n + \tau)) &= f(u_{i,j}^n) + \tau \partial_u f(u_{i,j}^n) \partial_t u_{i,j}^n + \mathcal{O}(\tau^2)
    \\
    &\approx f(u_{i,j}^n) + \tau \partial_u f(u_{i,j}^n) \left[ \kappa \Delta u_{i,j}^n + f(u_{i,j}^n) \right]
    \end{align}
If we set $r_x:= \kappa {\tau}/{h_x^2}$ and $r_y:= \kappa {\tau}/{h_y^2}$, this leads in eqn. \eqref{eq:cranknic} using eqn. \eqref{eq:taylorf} to 
    \begin{align}
   \notag &  u_{i,j}^{n+1}   - \frac 1 2 r_x \left[ u_{i,j-1}^{n+1} - 2 u_{i,j}^{n+1} + u_{i,j+1}^{n+1}                          \right] 
                        - \frac 1 2 r_y \left[ u_{i-1,j}^{n+1} - 2 u_{i,j}^{n+1} + u_{i+1,j}^{n+1} \right]
    \\
    &= u_{i,j}^n + \frac 1 2 r_x \left[ u_{i,j-1}^n - 2 u_{i,j}^n + u_{i,j+1}^n                                            \right] 
                  + \frac 1 2 r_y \left[ u_{i-1,j}^n - 2 u_{i,j}^n + u_{i+1,j}^n \right]
  + \tau f(u_{i,j}^n) + \frac 1 2 \tau^2 \partial_u f(u_{i,j}^n) \left[ \kappa \Delta u_{i,j}^n + f(u_{i,j}^n) \right].
    \end{align}
Thus, $\frac 1 2 \tau^2 \partial_u f(u_{i,j}^n) \left[ \kappa \Delta u_{i,j}^n + f(u_{i,j}^n) \right]$ is negligible for small values of $\tau$, which leads to the system
\begin{align}\label{eq:cn}
       \notag &-\frac 1 2r_x u_{i,j-1}^{n+1}-\frac 1 2r_x u_{i,j+1}^{n+1}+(1 + r_x + r_y)u_{i,j}^{n+1} -\frac 1 2 r_y u_{i-1,j}^{n+1} -\frac 1 2r_y u_{i+1,j}^{n+1} 
        \\
        &= \frac 1 2r_x u_{i,j-1}^n+\frac 1 2r_x u_{i,j+1}^n+(1 - r_x - r_y)u_{i,j}^n +\frac 1 2 r_y u_{i-1,j}^n +\frac 1 2r_y u_{i+1,j}^n  + \tau f (u_{i,j}^n) \, .
\end{align}
The system \eqref{eq:cn} leads to a linear equation system
\begin{align}
    A q_{n+1} = B q_n + \tau f_n \qquad \text{with, e.g.,} q_{n+1} = \left( [u_{0,0}^{n+1},...,u_{l_y,0}^{n+1}],[u_{0,1}^{n+1},...,u_{l_y,1}^{n+1}],....,[u_{0,l_x}^{n+1},...,u_{l_y,l_x}^{n+1}] \right)^T.
\end{align}
The vectors $q_n$ and $f_n$ are defined analogously, where $l_x$ and $l_y$ indicate the number of discretization points in $x$ and $y$ direction with respect to $\Omega$. The square and non singular matrices $A$ and $B$ are defined to contain the Neumann boundary conditions, which are implemented by, e.g. $u_{k+1,j}^{n+1}=u_{k,j}^{n+1}$ if $u_{k,j}^{n+1}$ lies on the boundary of the domain $\partial \Omega$. The previous refers to the solution of the PDE system of the state variable $u$. The adjoint system must be solved backwards in time, which leads to the approach
\begin{equation}
    \frac{z_{i,j}^{n-1}-z_{i,j}^n}{-\tau}=-\frac 1 2 \left[ \kappa \Delta z_{i,j}^{n-1} + \kappa \Delta z_{i,j}^n + p(u_{i,j}^{n-1}) + p(u_{i,j}^n) \right]
\end{equation}
where the $p(u_{i,j}^{n})$ contains componentwise the corresponding discretized terms of $\partial_{u_j} g + \sum_{k=1}^m z_k \partial_{u_j} f_k $. Proceeding analogously as before yields the system
\begin{align}
     \notag   & -\frac 1 2r_x z_{i,j-1}^{n-1}-\frac 1 2r_x z_{i,j+1}^{n-1}+(1 + r_x + r_y)z_{i,j}^{n-1} -\frac 1 2 r_y z_{i-1,j}^{n-1} -\frac 1 2r_y z_{i+1,j}^{n-1} 
        \\
        &= \frac 1 2r_x z_{i,j-1}^n+\frac 1 2r_x z_{i,j+1}^n+(1 - r_x - r_y)z_{i,j}^n +\frac 1 2 r_y z_{i-1,j}^n +\frac 1 2r_y z_{i+1,j}^n  + \tau p (u_{i,j}^n) \, .
\end{align}
Thus, when solving the linear system of equations, the matrices $A$ and $B$ can also be used, due to the same Neumann boundary conditions. 

\subsection{Finite element method for the SEIR-model}

An alternative to the Crank-Nicholson method which is used in the adjoint method, we also present a version of the finite element method which produces similar results. By plugging eqn. \eqref{eq:f(u)_} in eqn. \eqref{eq:u_PDE_}, the PDE system reads as follows:
\begin{subequations}
\label{eq__PDE_new}
    \begin{alignat}{6}
    \partial_t u_1 &= \kappa\, \Delta_{x,y} u_1 -\beta(t)u_1u_3,                   & (x,y) \in \Omega;\qquad&        u_1(t=0) &= u_{1,0}; \qquad& \partial_\nu u_1 &= 0 \, , \qquad&(x,y) \in \partial \Omega \,  \\
    \partial_t u_2 &= \kappa\, \Delta_{x,y} u_2 +\beta(t)u_1u_3-\theta u_2, \qquad & (x,y) \in \Omega;\qquad&        u_2(t=0) &= u_{2,0}; \qquad& \partial_\nu u_2 &= 0 \, , \qquad&(x,y) \in \partial \Omega \,  \\
    \partial_t u_3 &= \kappa\, \Delta_{x,y} u_3 +\theta u_2 - \gamma u_3,          & (x,y) \in \Omega;\qquad&        u_3(t=0) &= u_{3,0}; \qquad& \partial_\nu u_3 &= 0 \, , \qquad&(x,y) \in \partial \Omega \,  
\end{alignat}
\end{subequations}
As the diffusion and ODE parts in eqs. \eqref{eq__PDE_new} are handled by two different schemes, we can make use of an \textit{operator splitting} (cf. MacNamara \cite{Mac16}) to solve the system. E.g., for one time step $\Delta t$, this procedure is as follows:
\begin{enumerate}
\item[(1)] Solve $\partial_t u_i=\kappa \,\Delta_{x,y} u_i$ for $t=\frac {\Delta t} 2$ with the corresponding initial and boundary conditions for $i=1,2,3$. 
\item[(2)] Solve $\partial_t u_i=f_i(u)$ for $t=\Delta t$ with the corresponding initial and boundary conditions for $i=1,2,3$. 
\item[(3)] Solve $\partial_t u_i=\kappa\, \Delta_{x,y} u_i$ for $t=\frac {\Delta t} 2$ with the corresponding initial and boundary conditions for $i=1,2,3$. 
\end{enumerate}
The equation in (2) is a simple ODE equation which can be solved by any standard solver, e.g. the Euler method or the method of Runge-Kutta. To solve the equation in (1) and (3) on the domain $\Omega$ with initial and homogeneous Neumann boundary conditions on $\partial \Omega$, we consider its weak form gained by multiplication with a test function $v\in H^1_0(\Omega)$; this means, that $v$ vanishes at the boundary, i.e., $v\equiv 0$ at $\partial \Omega$. The weak form reads as follows:Instead of $u\in \mathcal C^{2,1}$, we now aim to find a function $u_i\in H^1(\Omega)$ solving
\begin{align}
\notag a(u_i,v):&=\underset{\Omega}\int\,  u_i \, v\, d\omega + \underset{\Omega}\int\, \nabla u_i\, \nabla v \,d\omega -\underset{\partial\Omega}\int\, u_i \, v\, d\lambda\\\notag&=\underset{\Omega}\int\,  u_i \, v\, d\omega +\underset{\Omega}\int\, \nabla u_i\, \nabla v \,d\omega\\&=0.
\end{align}
This infinite-dimensional problem has to be solved numerically by discretization. We aim to find a solution $u_{i,h}$ in a finite-dimensional subspace $V_h$ solving 
\begin{equation}
a(u_{i,h},v_h)=a_1(u_{i,h},v_h)+a_2(u_{i,h},v_h)=\underset{\Omega}\int\, u_{i,h} \, v_h\, d\omega +\underset{\Omega}\int\, \nabla u_{i,h}\, \nabla v_h \,d\omega=0\,.
\end{equation}
We define the subspace $V_h$ on the chosen grid and linearly independent basis functions $\phi_{j}$ piecewise over subregions $\Omega_k=[x_1,x_2]\times[y_1,y_2]\subset \Omega$:
\begin{align}
V_h={\left\{u_h=\sum_k \sum_{j=1}^4 c_j^{(k)} \,\phi_j^{(k)}(x,y)  \right\}}
\end{align}
where
\begin{subequations}
\begin{align}
\phi_1^{(k)}(x,y)&=\frac{(x-x_2)(y-y_2)}{(x_1-x_2)(y_1-y_2)}\\
\phi_2^{(k)}(x,y)&=\frac{(x-x_2)(y-y_1)}{(x_1-x_2)(y_2-y_1)}\\
\phi_3^{(k)}(x,y)&=\frac{(x-x_1)(y-y_2)}{(x_2-x_1)(y_2-y_2)}\\
\phi_4^{(k)}(x,y)&=\frac{(x-x_1)(y-y_1)}{(x_2-x_1)(y_2-y_1)}
\end{align}
\end{subequations}
for $(x,y)\in\Omega_k$; otherwise, those functions vanish, i.e., $\phi_j^{(k)}(x,y)\equiv 0$ for $(x,y)\not \in \Omega_k$, $j=1,2,3,4$. Then the weak form $a(u_{i,h},v_h)=0$ reads as follows:
\begin{align}
a(u_{i,h},v_h)=a\left(\sum_k \sum_{j=1}^4 c_j^{(k)}\,\phi_j^{(k)}(x,y),\, \phi_{j^*}^{({k^*})}(x,y)\right)=0\\
\intertext{and, due to the linearity of $a$,}
\sum_k \sum_{j=1}^4 a\left(\phi_j^{(k)}(x,y),\,\phi_{j^*}^{({k^*})}(x,y)\right) c_j^{(k)}=0
\end{align}
Then the stiffness matrices $A$ and $B$ are defined by 
\begin{subequations}
\begin{alignat}{2}
A_{nm}&=a_1\left(\phi_j^{(k)}(x,y),\,\phi_{j^*}^{({k^*})}(x,y)\right)\\
B_{nm}&=a_2\left(\phi_j^{(k)}(x,y),\,\phi_{j^*}^{({k^*})}(x,y)\right)
\end{alignat}
\end{subequations}
where $n$ represents the row corresponding to $(j^*,k^*)$ and $m$ the column corresponding to $(j,k)$, which depends on the chosen order within the matrices. More information about this can e.g. be found in \cite{Eva99}. The linear equation system with a mass matrix
\begin{align}
A \,\partial_t u_i+ B\, u_i = 0
\end{align}
 can be solved by any scheme; e.g., a 4-step Runge-Kutta scheme (cf. e.g. \cite{Dor80}). 

\section{Optimization and numerical nethods}

\subsection{Discretization of the domain}

To define the domain $\Omega$ for the partial differential equation model, we consider the geographical data for the district (D), the association communities (AC) and the municipalities (M).
  Using the free online data filtering tool \textit{Overpass Turbo} \cite{opt} from \textit{OpenStreetMap}, we extracted the relevant geographical data and created the relevant matrices assuming the map segment as rectangular, which appears reasonable due to the very small size of the investigated window. The relevant data for the starting values of the locations are then equally distributed across all relevant subdomains of $\Omega$. The size of the window is $L_x \times L_y= 39.23$ km $\times$ $56.05$ km.
  Using a step size of $h_x=\frac {L_x}{100}$ and $h_y=\frac {L_y}{100}$, the discretization of the area yields the following $101x101$ matrices.

\begin{figure}[H]
\begin{center}

		\begin{subfigure}{\textwidth}
		\centering
		\includegraphics[width=.56\linewidth]{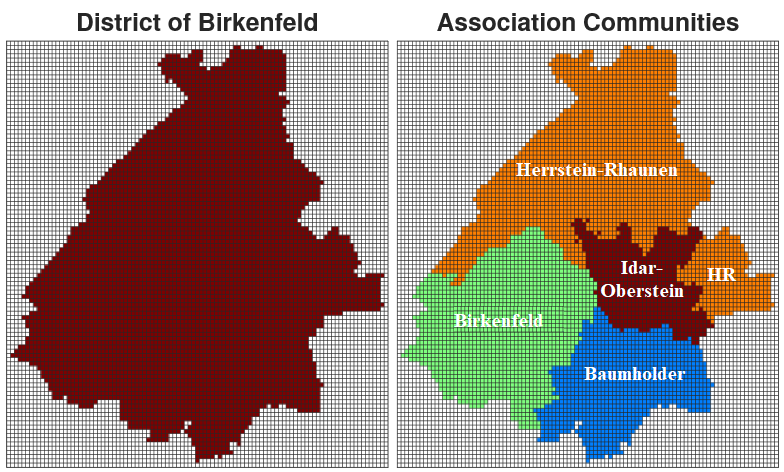}
		\end{subfigure}
		\centering
	\caption{Left: Discretizations of the whole district. Right: Discretization of the association communities.}
	\label{fig:data1B}
	\end{center}
\end{figure}

\begin{figure}[H]
\begin{center}
		\begin{subfigure}{.28\textwidth}
		\centering
		\includegraphics[width=1\linewidth]{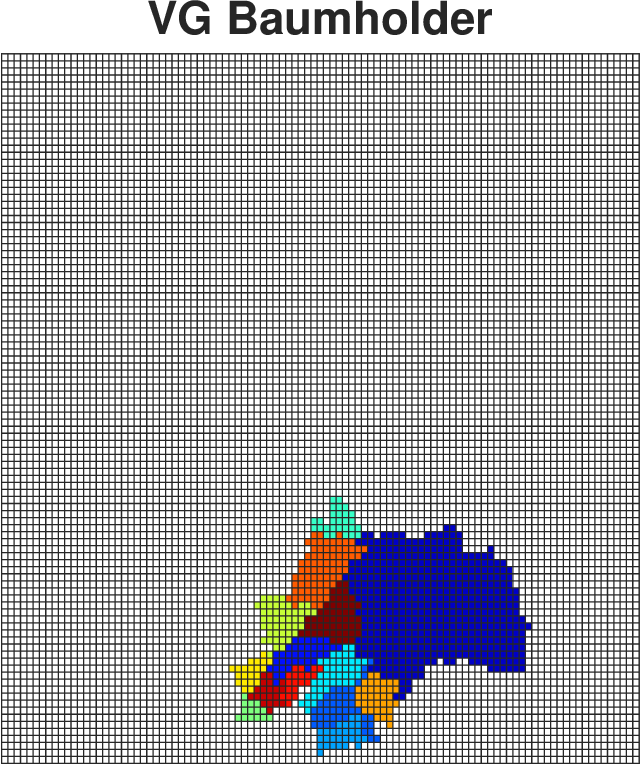}
		\end{subfigure}
		\begin{subfigure}{.28\textwidth}
		\centering
		\includegraphics[width=1\linewidth]{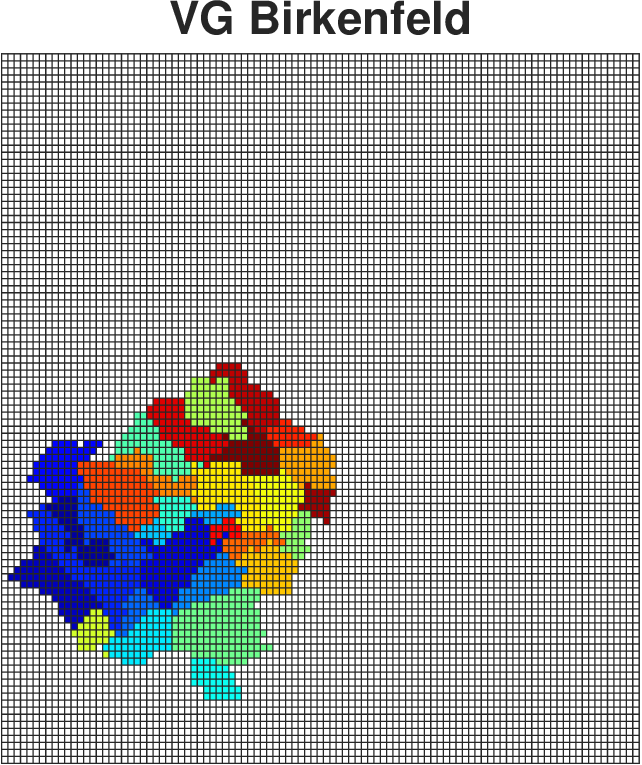}
		\end{subfigure}\vspace{0.15cm}\\
		\begin{subfigure}{.28\textwidth}
		\centering
		\includegraphics[width=1\linewidth]{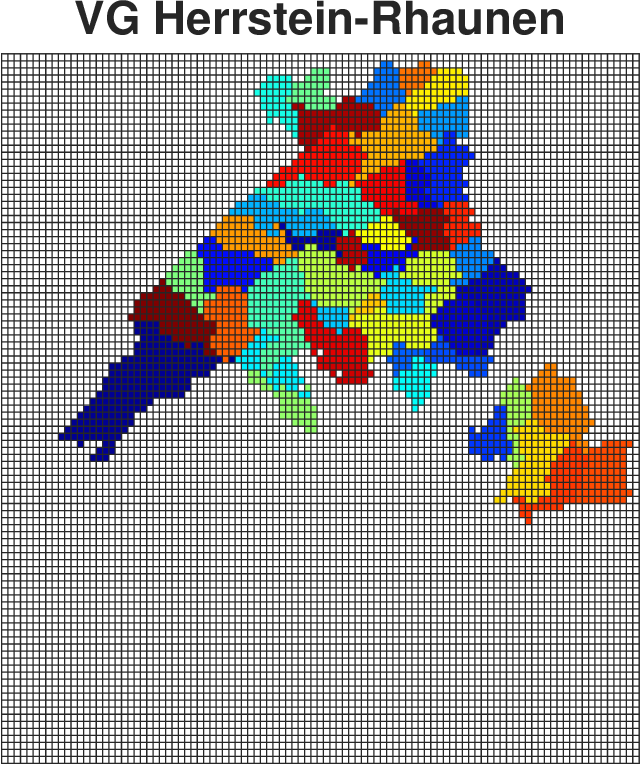}
		\end{subfigure}
		\begin{subfigure}{.28\textwidth}
		\centering
		\includegraphics[width=1\linewidth]{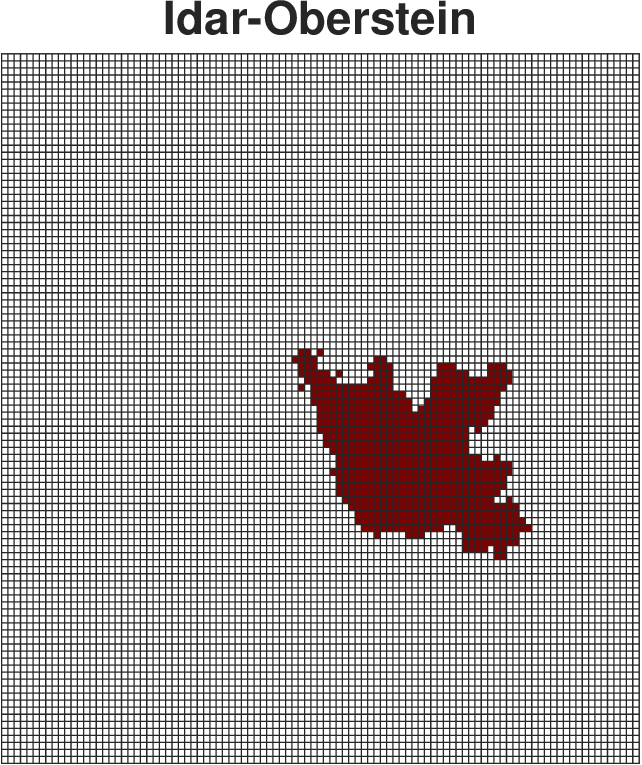}
		\end{subfigure}
	\caption{Discretization of all municipialities in the district of Birkenfeld; upper left: VG Baumholder, upper right: VG Birkenfeld, lower left: VG Herrstein-Rhaunen, lower right: city of Idar-Oberstein, which is classified as one municipiality.}
	\label{fig:data2}
			\end{center}
\end{figure}

\subsection{Target function}

In the following, we present the analysis for the reduced \textit{SEI}-model. The derivation for the other models is analogous. Furthermore, in order to avoid confusion, we rename $(u_S,u_E,u_I):=(u_1,u_2,u_3)$. The objective of this section is to present two methods for data fitting of the presented models to the dataset of Birkenfeld County in Germany. 

Therefore, we introduce an objective function $J : \mathbb{R}^d \times \mathcal{C}(\Omega,\mathbb{R}^m)^{2,1} $ defined by
\begin{equation}
\label{eq:J}
    J \left( \chi,u_0 \right) = \frac{w_0}{2} \lVert \delta \beta(t) u_S u_I - u_I^{data} \rVert_{L^2_V}^2 
    + \frac{w_1}{2} \lVert \chi-\tilde \chi \rVert_2^2 
    + \frac{w_2}{2} \sum_{j=1}^m \lVert u_{I,j,0} - u_{I,j,0}^\text{data} \rVert_{L^2_\Omega}^2 \, .
\end{equation}
On the one hand we make use of the $L^2$-norm with respect to $V$ and $\Omega$, defined for example by $\lVert g \rVert_{L^2_X} = \left( \int_X g^2 dx \right)^{1/2}$, on the other hand we use the Euclidean norm $\lVert x \rVert_2 = \left( 
\sum_{j=1}^d x_j^2 \right)^{1/2}$. 
The $L_V^2$-norm in $J$ involves fitting the model $u$ to the respective data sets, whereby $u_I^{data}: V \rightarrow \mathbb{R}^m$ is an interpolation through the data points. Interpolation is performed linearly with respect to the time axis for each grid point. We decided to fit the model to the data point by point in the $L_V^2$ norm. The reason for this is the necessary analysis in the adjoint method in subsection \ref{chap_adj}. 

The Euclidean norm in the objective function corresponds to a regularization term and contains the parameters of the model $u$, which have to be fitted to the data. For example, if we assume an unknown constant transmission rate $\beta(t) \equiv \beta$ and a diffusivity $\kappa$, we have $\chi=(\beta,\kappa,\delta)$. 
The $L_\Omega^2$-norms are to match the initial conditions of $u$ to the corresponding initial guess $u_{j,0}^{data}$, which is in the following simulations assumed to be $0$. The same applies to the initial guess of $\tilde{\chi}$.
In addition, the objective function includes weights $w_j$ with $j=0,1,2$, whose choice will be explained later.
The goal is now to minimize $J$ while satisfying the model constraints, i.e.
\begin{equation}
\label{min}
    \min_{\chi,u_0} \, J \, , \qquad \text{subject to PDE system } \eqref{eq:u_PDE_} \, .
\end{equation}
Note, that we only use informations concerning the daily reported new infected individuals. Thus, the data include the percentage of daily new infected individuals which shall be fitted to the incidence term $\beta(t) u_S u_I$ in the reduced and normalized \textit{SEI}-model. In addition, one must assume that only a fraction of those actually infected are reported. Therefore, we introduce a detection rate $\delta$, which is unknown.


\subsection{Models, parameter bounds and initial values}

Allowing the following respective constraints of the fitted parameters: \begin{align}[\beta_j,\kappa,u_0,\delta]\in[\R^+,[0,1],\R^4_+,[0,1]],\end{align}
the starting values at time $t_0$ for the detected cumulated infected $u^3_0$ can be taken from the statistics, while we assume no initial recovered person  $u_0^4=0$ and, for the initial amount of exposed persons, $u_0^2=u_0^3/2$ with similar reasoning as in \cite{HeS20}. The initial number of infected is then defined as $I_0=u_0^1+u_0^2+u_0^3+u_0^4$. 
\begin{table}[H]
    \centering
        \caption{Orders of magnitude of the initial values for adapting the model to the available data.}
    \label{tab___para}
    \begin{tabular}{cccccccc}
    param. &   $\beta_{j}$  & $\delta$ & $\kappa$ & $u_0^1$ & $u_0^2$ & $u_0^3$ & $u_0^4$
    \\
    \hline
    init. val. & $0.1$ & $0.5$ &$0.1$& $1/N_1$ &$1/N_2$& $1/N_3$&$3/N_4$ 
    \end{tabular}
\end{table}

\subsection{Metropolis algorithm}\label{chap_met}

The first presented method makes use of a Metropolis algorithm (cf. Metropolis et al. \cite{Met53}, Gelman et al. \cite{Gel96} or Gilks et al.\cite{Gil96}) for  estimation of parameters in the PDE system (\ref{eq__PDE_new}) according to the procedure described in Schäfer and Götz \cite{Sch20} and Heidrich, Schäfer et al. \cite{HeS20}. Using the parameter set $u_0$ as of Tab. \ref{tab___para} as starting conditions, we assign random draws $u_{new}$ from a normally distributed (and thus symmetric) proposal function $q$, i.e. $u_{new} \sim q(u_{new}|u_{i-1})$, in every iteration $i$. 
	
Using the previously defined $\hat J(u)$ as of eqn.\eqref{eq:J} as the target distribution, we calculate the approximative distribution by
\begin{equation}
\label{equa11}
\pi(u)=c \cdot \exp{\left(-{\frac{\hat J(\hat u)^2}{2 \sigma^2}}\right)},
\end{equation}
whereby $c$ is an arbitrary real value. For the acceptance probability, it follows 
\begin{equation}
\label{equa9}
p(\hat u_{new}|u_{i-1})=\min\left\{1, \frac{\pi(\hat u_{new})\cdot q(u_{i-1}|u_i)}{\pi(\hat u_{i})\cdot q(\hat u_{i}|\hat u_{i-1}))}\right\}=\min\left\{1, \frac{\pi(\hat u_{new})}{\pi(u_{i})}\right\}.
\end{equation}
In eqn. \eqref{equa9}, we can see that the value of $c$ is redundant, as it cancels out in the division. If the sample is accepted with the probability $p$, we set $\hat u_i=\hat u_{new}$; with the probability $1-p$, the sample is declined, meaning $\hat u=\hat u_{i-1}$ according to Rusatsi \cite{Rus15} or Schäfer and Götz \cite{Sch20}.

For parameter estimation using the Metropolis algorithm, we use algorithm \ref{algo3}.

	\begin{algorithm}[H]
\caption{Pseudocode for the Metropolis algorithm.}
\label{algo3}
    \begin{algorithmic}[1]
        \STATE $\pi,\hat u^{\text{data}} \gets \text{load initial values for $\pi$ and data}$
        \STATE $x,z \gets \text{solve PDE for state variable}$
        \STATE $\hat J \gets \text{compute objective function regarding $\pi$}$ 
        \STATE $\text{$\sigma$} \gets \text{standard distribution of the solution, i.e. $I$ over time}$
        \STATE $s \gets \text{set step size (standard deviation) for the algorithm, e.g. $s:=\pi/100$}$
        \REPEAT
            \STATE $\pi_{old} \gets \text{$\pi$ from previous draw}$
            \STATE $ \pi_{new} \gets \pi \sim \mathcal N (\pi_{old},s)$
            \STATE $x, z, J(\hat \pi_{new}) \gets \text{update depending on} \; \pi$ 
            \STATE $\alpha \gets \min\left\{1, \exp\left({{\hat  J(\pi_{old})^2-\hat J(\pi_{new})^2}/{2 \sigma^2}}\right) \right\}$
            \STATE $\pi_{new} \gets \text{$\hat \pi_{new}$ with probability $\alpha$ and $\pi_{new}:=\pi$ with probability $1-\alpha$}$
        \UNTIL \text{maximum value of draws is reached}
        \STATE $\pi^*,x^*,\hat J^* \gets \text{means of all $\pi,x,J$}$
    \end{algorithmic}
\end{algorithm}

\subsection{Parameter estimation via adjoint functions}\label{chap_adj}

For parameter estimation via adjoint functions $z \in \mathcal{C}(V,\mathbb{R}^m)^{2,1}$ we use them in conjunction with a Lagrange function defined as $\mathcal{L}: \mathbb{R}^d \times \mathcal{C}(\Omega,\mathbb{R}^m)^{2,1} \times \mathcal{C}(V,\mathbb{R}^m)^{2,1} \times \mathcal{C}(V,\mathbb{R}^m)^{2,1} \rightarrow \mathbb{R}$, fulfilling
\begin{equation}
\label{eq:lagrange}
    \mathcal{L}(\chi,u_0,u,z) = 
                                J(\chi,u_0)
                                + \sum_{j=1}^m \int_V z_j \left( f_j(u)+\kappa \Delta_{x,y,} u_j - \partial_t u_j \right) \, d\omega dt \, .
\end{equation}
At a possible minimum $(\chi_*,u_{0*},u_*,z_*)$ of problem \eqref{min} must apply
\begin{equation}
    0 = \nabla \mathcal{L} = \left(\partial_{\chi_1} \mathcal{L}, ... , 
    \partial_{\chi_d} \mathcal{L}, 
    \partial_{u_{1,0}} \mathcal{L} , \dots, \partial_{u_{m,0}} \mathcal{L},
    \partial_{u_1} \mathcal{L} , \dots , \partial_{u_m} \mathcal{L} ,
    \partial_{z_1} \mathcal{L} , \dots ,  \partial_{z_m} \mathcal{L}  \right) \, .
\end{equation}
The derivatives in directions representing functions are determined with the help of G\^ateaux derivatives. The application of the optimal control theory and Pontryagin's maximum (minimum) principle provide the following optimality conditions:
\begin{tabular}{rlll}
    \\
     (i)    &   $0 = \partial_{\chi_k} \mathcal{L}$,     
            &   $k=1,\dots,d$,
            &   (Scalar Optimality Condition)  
     \\
     \\
     (ii)   &   $u_{j,0}=u_{j,0}^{data}-\frac{z_j(x,y,0)}{w_2}$,
            &   $j = 1,\dots,m$,
            &   (Optimal Initial Conditions)  
     \\
     \\
     (iii)  &   $\partial_t z_j = - \left( \partial_{u_j} g + \sum_{k=1}^m z_k \partial_{u_j} f_k            + \kappa \Delta_{x,y} z_j \right)$,
            &   $g:=\frac{w_0}{2} \left( \delta \beta(t) u_S u_I - u_I^{data} \right)^2 $
            &   (Adjoint Equations)
            \\
            &   $z=0$, 
            &   $t=t_{end}$,
            &   (Transversality Conditions)
            \\
            &   $\partial_\nu z = 0$, 
            &   $(x,y) \in \partial \Omega$,
            &   (Neumann Boundary Conditions) .
            \\
            \\
\end{tabular}

It should be noted that the present conditions are given as $g$ for the reduced \textit{SEI}-model. The adjoint equations PDE system has be solved numerically backward in time. Let us assume a time--dependent transmission rate of the form \eqref{eq:beta}.
This leads to $\chi=\left( \beta_0,\beta_1,\beta_2,\kappa,\delta \right)$ and to a corresponding gradient
\begin{subequations}
\label{eq_grad_}
\begin{alignat}{2}
    \partial_{\beta_k} \mathcal{L} &=
    w_1 \left( \beta_k - \tilde{\beta}_k \right) + \int_{\mathcal{I}_k} \int_\Omega  w_0 \delta u_S u_I \left( \delta \beta(t) u_S u_I - u_I^{data} \right) + \left( z_2 - z_1 \right) u_S u_I \, d\omega \,dt \, , 
    \\
    & \text{where $\mathcal{I}_0 = [0,t_0]$, $\mathcal{I}_1 = [t_0,t_1]$ and  $\mathcal{I}_2 = [t_1,t_{end}]$ for $k=0,1,2$} \, ,
    \\
    \partial_\kappa \mathcal{L} &= w_1 \left( \kappa - \tilde{\kappa} \right) 
    + \sum_{j=1}^3 \int_V z_j \Delta_{x,y} u_j \, d\omega dt \, , 
    &&
    \\
    \partial_\delta \mathcal{L} &=  w_1 \left( \delta - \tilde{\delta} \right)
    + w_0 \int_V \beta(t) u_S u_I \left( \delta \beta(t) u_S u_I - u_I^{data} \right) \, d\omega \,dt \, . &&
\end{alignat}
\end{subequations}
As adjoint equations we get in this case
\begin{subequations}
\label{eq:adjoint}
\begin{align}
    \partial_t z_1 &= - \left( w_0 \delta \beta(t) u_I \left( \delta \beta(t) u_S u_I - u_I^{data} \right)
    + \left( z_2 - z_1 \right) \beta(t) u_I + \kappa \Delta_{x,y} z_1 \right) \, , 
    \\
    \partial_t z_2 &= -\left( \theta \left( z_3 - z_2 \right) + \kappa \Delta_{x,y} z_2 \right) \, , 
    \\
    \partial_t z_3 &= -\left( w_0 \delta \beta(t) u_S \left( \delta \beta(t) u_S u_I - u_I^{data} \right) + \left( z_2- z_1 \right) \beta(t) u_S - \gamma z_3  + \kappa \Delta_{x,y} z_3 \right) \, .
\end{align}
\end{subequations}

For parameter estimation, we use algorithm \ref{alg:adj}. 
\begin{algorithm}[H]
	\caption{Pseudocode for the parameter estimation via adjoint functions.}
	\label{alg:adj}
    \begin{algorithmic}[1]
        \STATE $\beta_k,\kappa,\delta,u_I^{DATA},u_0^{data} \gets \text{load initial values and data}$
        \STATE $u,z \gets \text{solve PDE for state variable and adjoint function}$
        \STATE $J, \; \nabla J \gets \text{compute objective function and gradient regarding $\chi=(\beta_0,...,\beta_k,\kappa,\delta)$}$ 
        \STATE $s_1 \gets \text{compute search direction regarding $\chi$ } \text{(Quasi-Newton (BFGS))}$
        \STATE $s_2 \gets (\tilde{u}_0 - u_0)$ compute search direction for $u_0$ with $\tilde{u}_{j,0}=u_{j,0}^{\textit{data}}-\frac{z_j(x,y,0)}{w_2}$
        \REPEAT
            \STATE $J_{old} \gets J$
            \STATE $\alpha \gets 1$
            \STATE $\chi \gets \chi+\alpha s_1$
            \STATE $u_0 \gets u_0 + \alpha s_2$
            \STATE $u,J \gets \text{update}$
            \REPEAT
            \STATE $\alpha \gets 0.5 \alpha$
            \STATE $\chi \gets \chi+\alpha s_1$
            \STATE $u_0 \gets u_0 + \alpha s_2$
            \STATE $u,J \gets \text{update}$
            \UNTIL $J \leq J_{old} + 0.001 \alpha s^T \nabla J_{old}$ \; (Armijo Rule)
            \STATE $z, \nabla J, s_1, s_2 \gets \text{update}$ 
        \UNTIL $\frac{\Vert J - J_{old} \Vert_2}{\Vert J_{old} \Vert_2} < \textsc{TOL}$
    \end{algorithmic}
\end{algorithm}
Regarding the optimization of $\chi$, a Quasi-Newton Broyden-Fletcher-Goldfarb-Shanno (BFGS) search direction is used. The initial condition $u_0$ is updated with a convex combination between old value and current 'optimal' value. The step size is mediated by the Armijo stepsize rule. In each optimization step, the PDE system for the state $u$ and adjoint $z$ variable must be solved. This is done by the Crank--Nicholson method presented above. The fact that the state variable must be solved forward and the adjoint variable backward in time also leads to the term "forward--backward sweep method".

\section{Numerical results}

\subsection{Without penalty term (Metropolis)}

\begin{figure}[H]
\begin{center}
		\begin{subfigure}{.45\textwidth}
		\centering
		\includegraphics[width=1\linewidth]{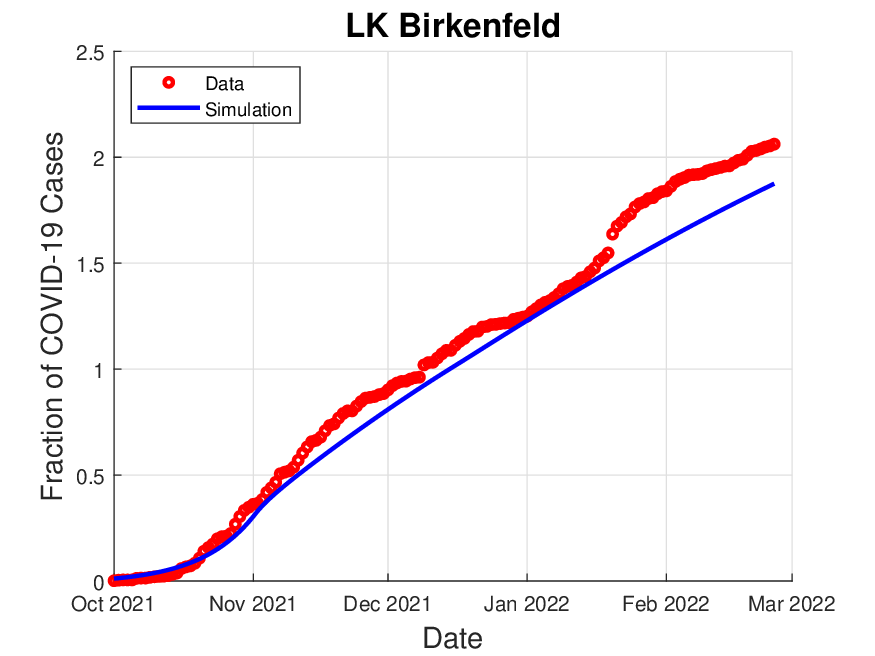}
		\end{subfigure}
		\begin{subfigure}{.45\textwidth}
		\centering
		\includegraphics[width=1\linewidth]{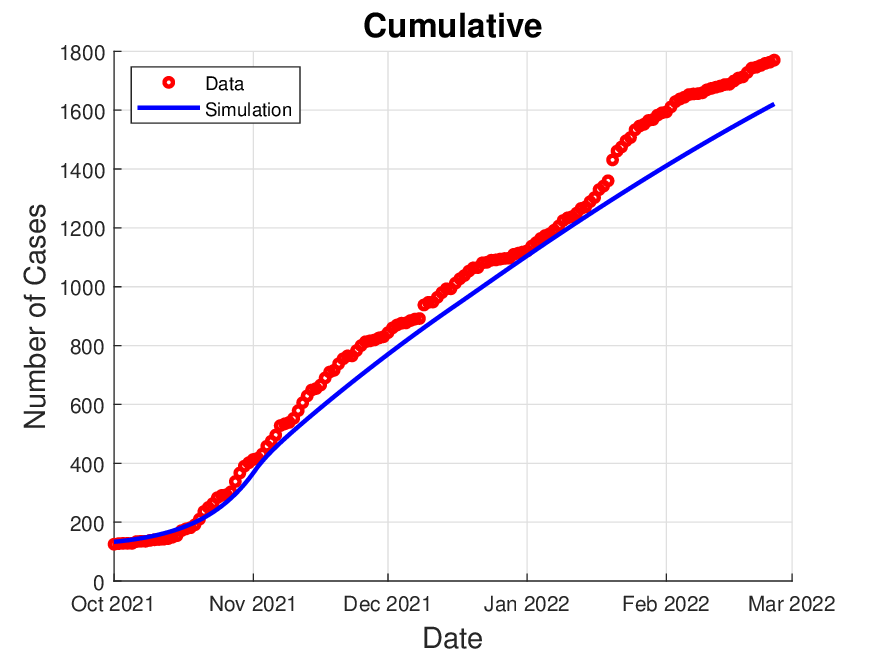}
		\end{subfigure}
	\caption{Result of the optimization with the Metropolis algorithm for the district of Birkenfeld with $w_1=w_2=0$.}
	\label{fig:data6}
	\end{center}
\end{figure}

\begin{figure}[H]
\begin{center}
		\begin{subfigure}{.45\textwidth}
		\centering
		\includegraphics[width=1\linewidth]{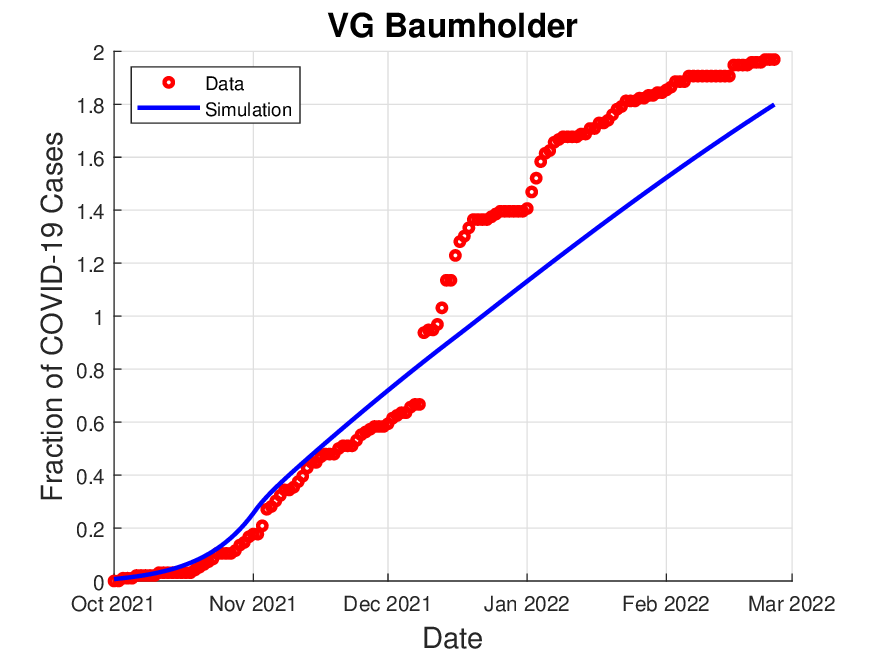}
		\end{subfigure}
		\begin{subfigure}{.45\textwidth}
		\centering
		\includegraphics[width=1\linewidth]{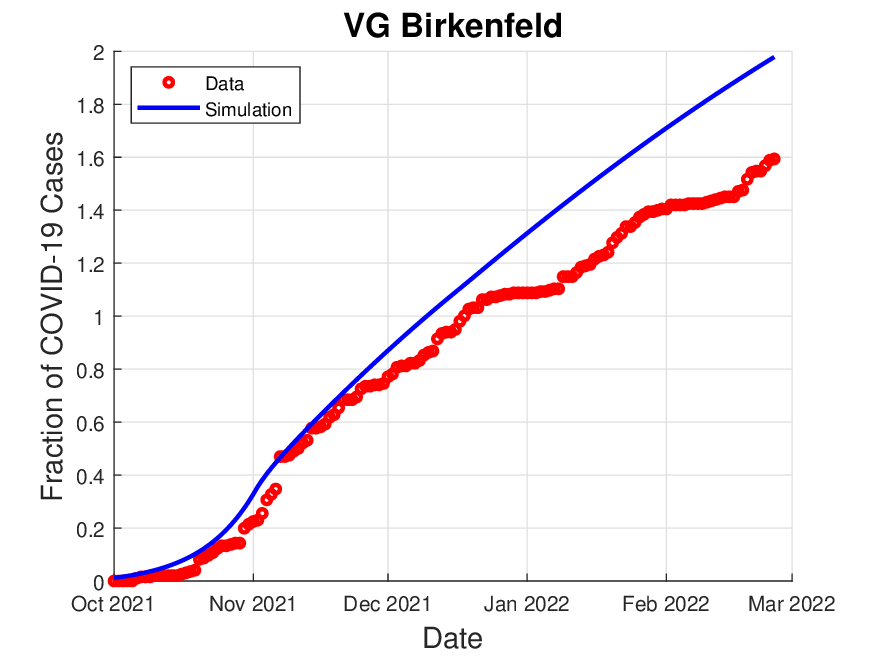}
		\end{subfigure}\\
		\begin{subfigure}{.45\textwidth}
		\centering
		\includegraphics[width=1\linewidth]{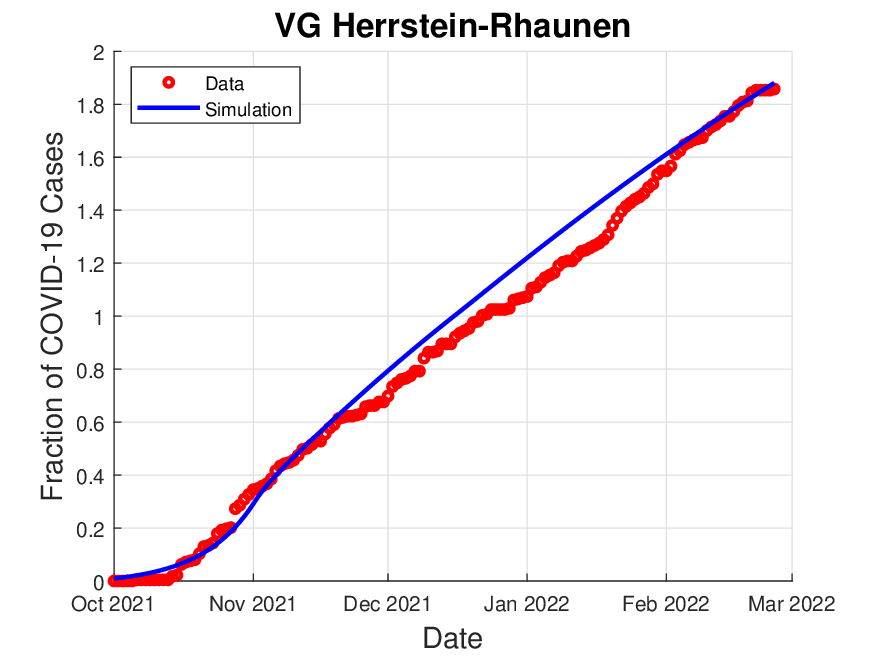}
		\end{subfigure}
		\begin{subfigure}{.45\textwidth}
		\centering
		\includegraphics[width=1\linewidth]{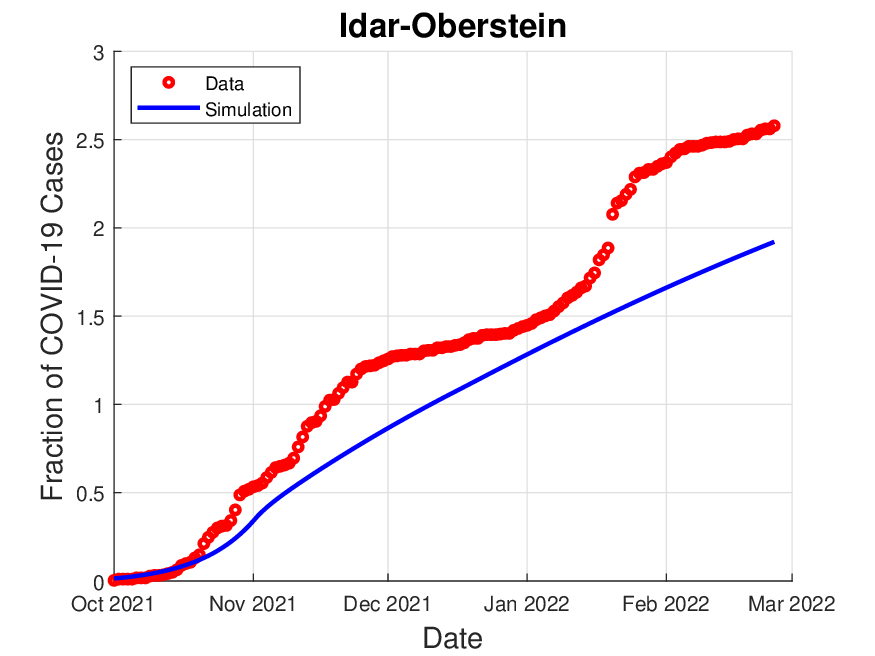}
		\end{subfigure}
	\caption{Result of the optimization with the Metropolis algorithm for the lower level administrative units (Verbandsgemeinden) with $w_1=w_2=0$.}
	\label{fig:data7}
			\end{center}
\end{figure}

Using the Metropolis algorithm as of chapter \ref{chap_met}, we firstly set $w_0=1$ and $w_1=w_2=0$ in eqn. \eqref{eq:J}.  In Fig. \ref{fig:data6}, the results for the district of Birkenfeld and in Fig. \ref{fig:data7}, the results for the four ACs are presented; the red dots represent the respective data, the blue line the outcome of the model. 

\subsection{With penalty term}

\begin{figure}[H]
\begin{center}
		\begin{subfigure}{.45\textwidth}
		\centering
		\includegraphics[width=1\linewidth]{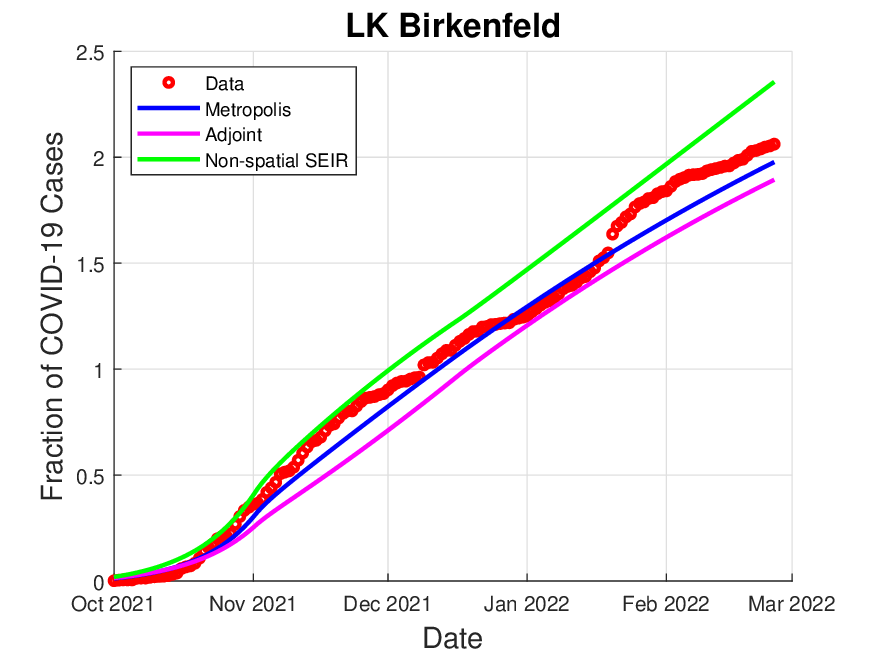}
		\end{subfigure}
		\begin{subfigure}{.45\textwidth}
		\centering
		\includegraphics[width=1\linewidth]{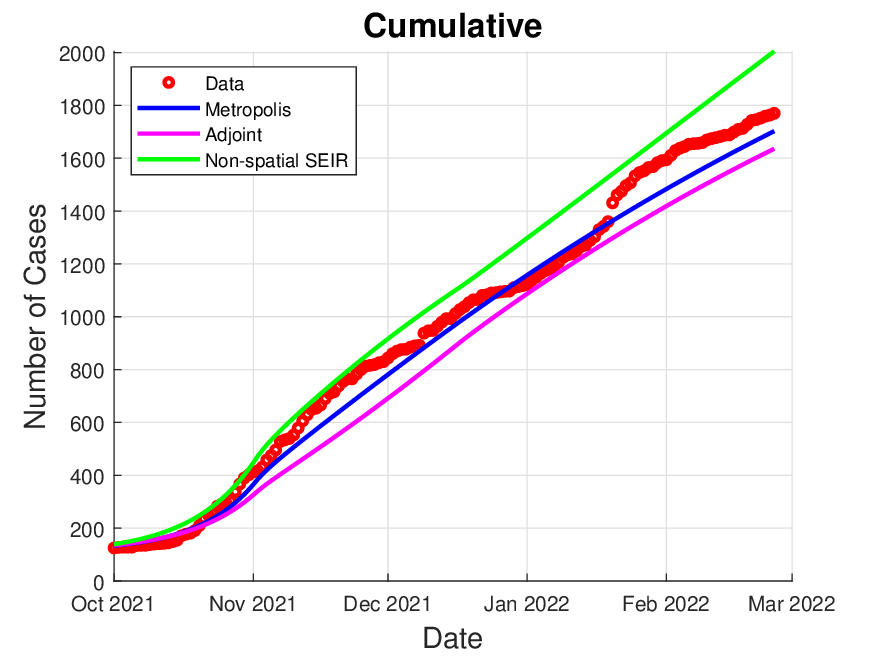}
		\end{subfigure}
	\caption{Result of the optimization with the various methods for the district of Birkenfeld with $w_1=w_2=10^{-5}$.}
	\label{fig:data8}
	\end{center}
\end{figure}

\begin{figure}[H]
\begin{center}
		\begin{subfigure}{.45\textwidth}
		\centering
		\includegraphics[width=1\linewidth]{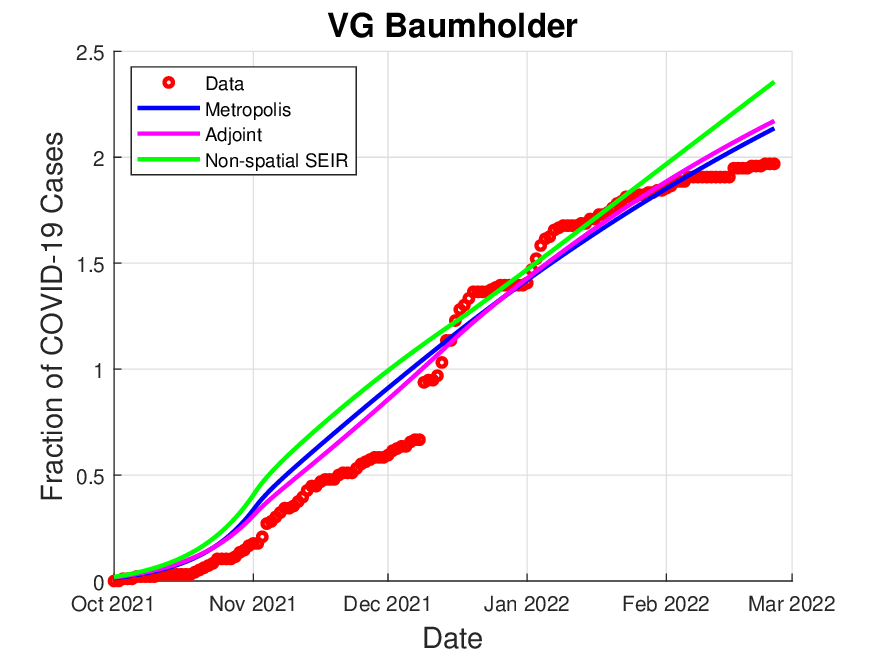}
		\end{subfigure}
		\begin{subfigure}{.45\textwidth}
		\centering
		\includegraphics[width=1\linewidth]{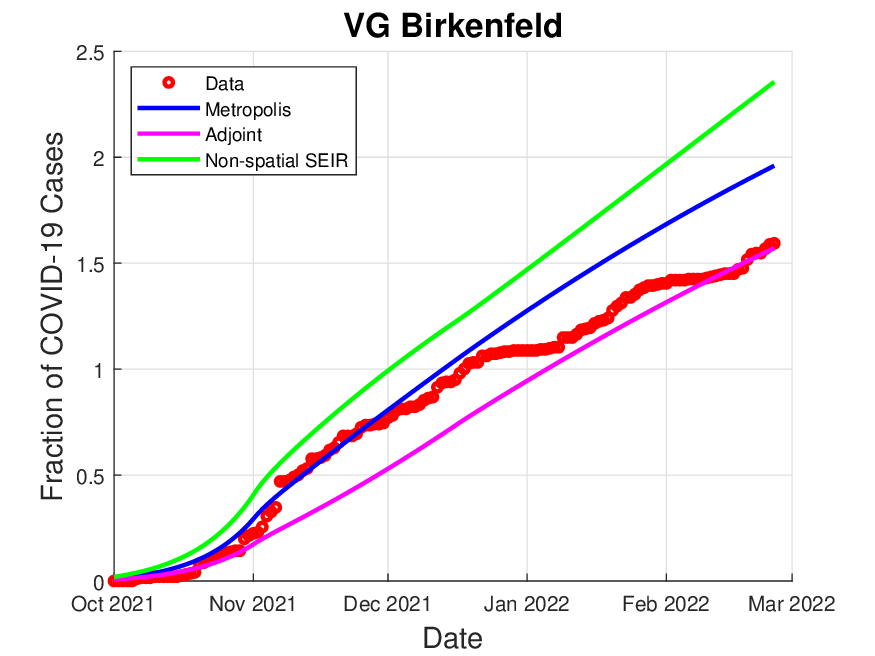}
		\end{subfigure}\\
		\begin{subfigure}{.45\textwidth}
		\centering
		\includegraphics[width=1\linewidth]{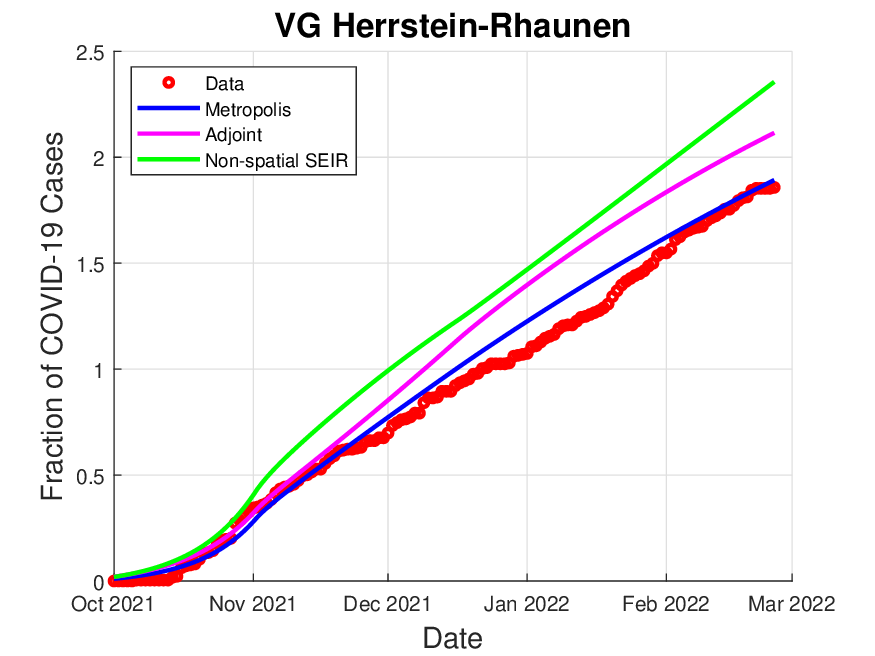}
		\end{subfigure}
		\begin{subfigure}{.45\textwidth}
		\centering
		\includegraphics[width=1\linewidth]{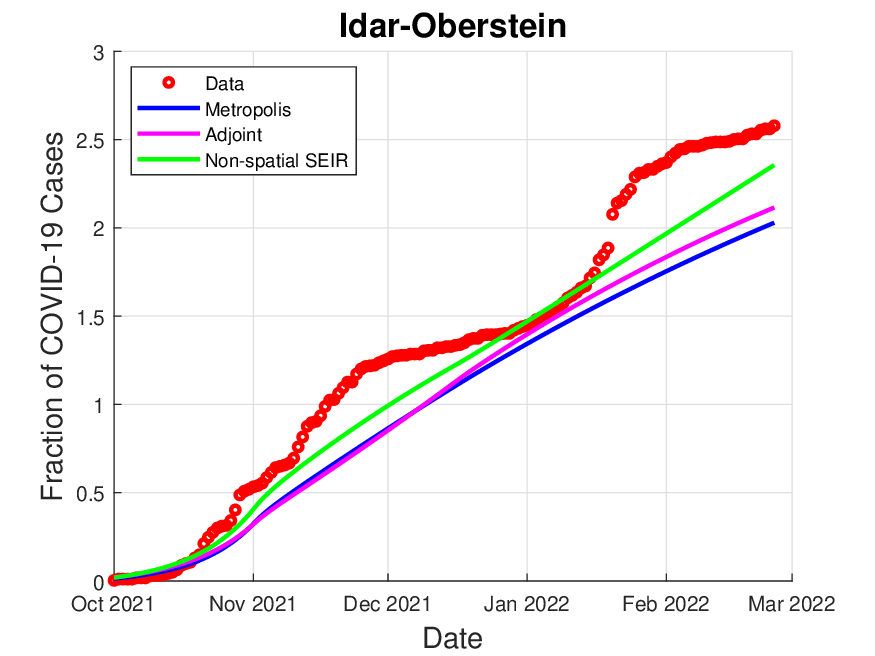}
		\end{subfigure}
	\caption{Result of the optimization with the various methodes for the lower level administrative units (Verbandsgemeinden) with $w_1=w_2=10^{-5}$.}
	\label{fig:data9a}
			\end{center}
\end{figure}

We now include a penalty term; in eqn. \eqref{eq:J} we therefore choose $w_0= 1$ and $w_1=w_2=1\cdot10^{-5}$, which guarantees a {convex} problem, cf. Heidrich, Schäfer et al. \cite{HeS20}. In Fig. \ref{fig:data8}, the results for the district of Birkenfeld and in Fig. \ref{fig:data9a}, the results for the four associated communities are presented; the red dots represent the respective data, the blue line the outcome of the model. We also compare the results to those of a standard \textit{SEIR}-model with the usual parameter estimation (using the Metropolis algorithm).

\subsection{Comparison}

In the results of the parameter estimation as shown in Tab. \ref{tab_num_res4}, it is notable that across all simulations, we find $\beta_1 \approx \beta_2$, which indicates that the more severe lockdown restrictions from December 17 onwards might be overlaid by the rising number of festivities during Christmas. The adjustment regarding the detection rate shows that, according to the model, the  actual  number  of  infected  people  is  approximately  three to five  times  higher  than registered, which is in line with previous findings.

The parameter values for the transmission parameters $\beta_i$ are relatively consistent across the different methods and parameter settings. The values $\beta_0>\beta_1\approx \beta_2$ show that the so-called 'light' lockdown in November 2020 had the most significant effect on the case values, while the more 'severe' lockdown before Christmas did not have a significant effect. A reason for this can be Christmas itself during which the amount of contacts has rised by nature, it also has to be noted that the measurement of case numbers during or after Christmas was not consistent, which might lead to delays in the infection data. 

The parameter $\kappa$ is ranging significantly between 0.05 and 0.12, but correlations to the higher initial values of the infected, as well as the detection rate $\delta$, vary across the different methods and parameter settings, depending not only on the weights $w_1$ and $w_2$ but also the chosen $\kappa$. The results indicate that there is some potential variation in the optimal parameter values and parameter cross-correlation of parameters, e.g. of the initial infection values and the detection rate, as well as to some extend in the diffusivity. This can also be seen in the target function $J(u)$ which, despite the deviations, shows only minor variations among the different methods and parameter settings.
    
The estimations for some of the districts is generally quite accurate, e.g. for Herrstein-Rhaunen and Baumholder. However, for the city of Idar-Oberstein, the area with most infections, the model underestimates the infection cases, which can be explained by the diffusivity of the model that neglects the urban structures. On the other side, the model overestimate the infection cases in Birkenfeld. In total, the amount of infected is slightly underestimated towards the end, while the Metropolis algorithm with $c=10^{-5}$ provides the most accurate results for the cumulative data. Compared to the \textit{SEIR}-model without diffusion, we see that both Metropolis and adjoint methods provide better target function values, as the \textit{SEIR}-model generally overestimates the infections in all districts except the largest (Idar-Oberstein).
\begin{table}[H]
    \centering
       \caption{Results for the Metropolis algorithm and the adjoint method, compared to an \textit{SEIR}-model without diffusion. For the Metropolis simulations, the standard deviation is given as an addition.}
   \label{tab_num_res4}
  \begin{tabular}{c|c|c|c|c}  
      &  Metropolis & Metropolis &  Adjoint& \textit{SEIR} \\
    $w_1=w_2$   &   0 & $10^{-5}$ &  $10^{-5}$& $10^{-5}$ \\
    \hline
    $\beta_0$   & $({0.228} \pm {0.005})$ d$^{-1}$& $({0.218} \pm {0.011}) $ d$^{-1}$&0.202 d$^{-1}$ &     0.206
   \\
     $\beta_1$  & $({0.097}\pm {0.003}) $ d$^{-1}$& $({0.099} \pm {0.003}) $ d$^{-1}$&0.109 d$^{-1}$ &   0.092
     \\
   $\beta_2$   &$({0.103} \pm {0.029} )$ d$^{-1}$ &$({0.097} \pm {0.027}) $ d$^{-1}$& 0.097 d$^{-1}$&    0.105
   \\
   $\kappa$     & $0.119\pm 0.004$  & $0.100\pm 0.005$ &0.102 & 0.000
   \\
   $I_{0}^{BA} $   & $3.393\pm0.077$ & $2.629\pm 0.101$&4.007&   4.139
   \\
   $I_{0}^{BI} $   & $11.570\pm 1.183$& $4.424\pm 0.258$&3.275&8.441
   \\
   $I_{0}^{HR} $   & $10.176 \pm 0.276$& $4.753\pm 0.2615$&6.350&
9.632
   \\
   $I_{0}^{IO} $   & $20.681\pm 0.682$& $8.024\pm 0.291$&14.630&12.207
   \\
   $\delta $   & $0.202\pm 0.005$ & $0.495\pm 0.010$ &0.397&0.444
   \\   \hline
   $J(u) $   & 0.4829&  0.4836 &0.4850 & 0.4858 
   \end{tabular}
\end{table}

\section{Discussion}

This study presented a reaction-diffusion model used to simulate the spatial spread of COVID-19, making use data down to the municipality level for parameter estimation. The \textit{SEIR}-model is based on a set of PDEs that describe the dynamics of susceptible, exposed, infected, and recovered individuals, as well as an incidence term that represents the transmission of the disease. The optimal control is based on a least square fit between the model output and the reported daily data. Two different approaches for the estimation of parameters and approximation of the infection data -- the Metropolis algorithm and the adjoint method -- were described and implemented, and their results were plotted and compared. 

Regarding the graphical and numerical results, all routines have provided  meaningful results. 
The models depict the infection values quite accurately in several subdistricts, yet slightly over- or underestimate them in others, which can partially be explained by non-homogeneous behaviour of cities compared to rural areas. On a local level, the quality of the estimations decreases -- as expected, as it cannot be assumed that a global model will apply perfectly to the behavior of single villages, especially as some of them had no or less than a handful of detected infected in the observed time interval. The parameter values, including initial infection values, detection rate, and diffusivity, vary across different methods and parameter settings, while the transmission-related parameters remain relatively consistent, and the target function are very similar. Compared to a non-spatial \textit{SEIR}-model, both Metropolis and adjoint method provide better results with respect to the target function $J(u)$.

It is important to note that the PDE model still provides valuable insights and information and manages to describe the diffusion in the epidemiological situation (very few and locally condensed initial cases in a more or less completely susceptible population). Further refinement of the model or the use of additional data sources may improve its accuracy in the future. Nevertheless, the accuracy of the model on the district and ACs level make it a valuable tool for understanding the spread of infections. Not only new variants of COVID-19, but also the possibility of future pandemics underscore the need for accurate modeling on local level.

\section*{Acknowledgements}
We would like to thank Karsten Schultheiß and his staff at Department 1 for Central Tasks and Issues of the Birkenfeld County. The department, among other things, is responsible for press and public relations and provided us with the data used.

\bibliography{main}

\end{document}